\documentclass[preprint,12pt]{article}

\newcommand{\ignore}[1]{}
\ignore{The text written between the parentheses here is ignored by latex}

\makeatletter
\def\ps@pprintTitle{%
 \let\@oddhead\@empty
 \let\@evenhead\@empty
 \def\@oddfoot{\reset@font\hfil\thepage\hfil}
 \let\@evenfoot\@oddfoot
}
\makeatother
\usepackage{amssymb}

\usepackage{amsmath}
\usepackage{algorithm}
\usepackage{algorithmic}
\usepackage{graphicx}
\usepackage{multirow}
\usepackage{lscape}
\usepackage{rotating}
\usepackage{latexsym,amsmath,amssymb,amsthm,bm,booktabs,authblk,multicol}
\usepackage{setspace}
\usepackage{epsf}
\usepackage{url}
\usepackage{color}

\def\@begintheorem#1#2{\par\bgroup{\sc #1\ #2. }\it\ignorespaces}
\def\@opargbegintheorem#1#2#3{\par\bgroup{\sc #1\ #2\ (#3). } \it\ignorespaces}
\def\@endtheorem{\egroup}
\newtheorem{theorem}{Theorem}[section]

\newtheorem{example}[theorem]{Example}
\newtheorem{proposition}[theorem]{Proposition}
\newtheorem{property}[theorem]{Property}

\newtheorem{conjecture}[theorem]{Conjecture}

\def \conv {{\rm conv}}

\def\ZZ{\mathbb{Z}}

\def\RR{\mathbb{R}}

\begin{document}

\title{\bf Primitive Zonotopes}
\date{ }

\author[a]{Antoine Deza}
\author[b]{George Manoussakis} 
\author[c]{Shmuel Onn}

\affil[a]{{\small 
McMaster University, Hamilton, Ontario, Canada}}
\affil[b]{{\small 
Universit\'e de Paris Sud, Orsay, France}}
\affil[c]{{\small Technion - Israel Institute of Technology, Haifa, Israel}}

\maketitle
\begin{abstract}
\noindent
We introduce and study a family of polytopes which can be seen as a generalization of the permutahedron of type $B_d$. We highlight connections with the largest possible diameter of the convex hull of a set of points in dimension $d$ whose coordinates are integers between $0$ and $k$, and with the computational complexity of multicriteria matroid optimization.
\end{abstract}


\section{Introduction}
We introduce and study lattice polytopes generated by the primitive vectors of bounded norm. These {\em primitive  zonotopes} can be seen as a generalization of the permutahedron of type $B_d$. We note that,  besides a large symmetry group,  primitive  zonotopes have a large diameter and many vertices relative to their grid size embedding. The article is structured as follows. In Section~\ref{eulerZ}, we introduce the primitive  zonotopes and some of their properties.  In Section~\ref{deltaZ} we derive lower bounds for the diameter of lattice polytopes and in Section~\ref{vZ} we determine the computational complexity of multicriteria matroid optimization.

 Finding a good bound on the maximal edge-diameter of a polytope in terms of its dimension and the number of its facets is not only a natural question of discrete geometry, but also historically closely connected with the theory of the simplex method.
Recent results dealing with the combinatorial, geometric, and algorithmic aspects of linear optimization include Santos' counterexample to the Hirsch conjecture, and Allamigeon, Benchimol, Gaubert, and Joswig's counterexample to a continuous analogue of the polynomial Hirsch conjecture. Kalai and Kleitman's upper bound for the diameter of polytopes was strengthened by Todd, and then by Sukegawa. Kleinschmidt and Onn's upper bound for the diameter of lattice polytopes was strengthened by Del Pia and Michini, and then by Deza and Pournin. For more details and additional results such as the validation that transportation polytopes satisfy the Hirsch bound, see~\cite{ABGM14,BSEHN14,BDF16,DM16,DP16,KK92,S12,S16,T14} and the references therein.

Multicriteria matroid optimization is a generalization of standard linear matroid optimization where each basis is evaluated according to several, rather than one criteria, and these values are traded-in by a convex function, see~\cite{MO,Onn,OR} and the references therein. In turns out that multicriteria matroid optimization can be reduced to solve several linear counterparts. In Section~\ref{vZ}, the largest number of such counterparts is shown to be precisely the number of vertices of some primitive zonotopes.

\section{Primitive zonotopes}\label{eulerZ}
\subsection{Definition}
The convex hull of integer-valued points is called a lattice polytope and, if all the vertices are drawn from $\{0,1,\dots,k\}^d$,  is refereed to as a lattice $(d,k)$-polytope.  For simplicity, we only consider full dimensional lattice $(d,k)$-polytopes. Given a finite  set $G$ of vectors, also called the generators, the zonotope generated by $G$ is the convex hull of all signed sums of the elements of $G$.  Searching for lattice polytopes with a large diameter for a given $k$, natural candidates include zonotopes generated by short integer vectors in order to keep the grid embedding size relatively small. In addition,  we restrict  to integer vectors which are pairwise linearly independent  in order to maximize the diameter. Thus, for $q=\infty$ or a positive integer, and $d,p$ positive integers, we consider the  {\em primitive  zonotope} $Z_q(d,p)$ defined as the zonotope  generated by the primitive integer vectors of $q$-norm at most $p$:
\begin{eqnarray*}
Z_q(d,p)
& = & \sum [-1,1]\{v\in\ZZ^d\,:\,\|v\|_q\leq p\,,\ \gcd(v)=1\,,\ v\succ 0\} \\
& = &\conv\left(\sum\{\lambda_v v\,:\, v\in\ZZ^d\,,
\|v\|_q\leq p\,,\ \gcd(v)=1\,,\ v\succ 0\}\ :\ \lambda_v=\pm1\right)\ .
\end{eqnarray*}
where $\gcd(v)$ is the largest integer dividing all entries of $v$, and $\succ$  the lexicographic order on $\RR^d$, i.e.  $v\succ 0$ if the first nonzero coordinate of $v$ is positive. Similarly, we consider $H_q(d,p)$ which is, up to  translation, the image of $Z_q(d,p)$ by a homothety of  factor $1/2$:
$$H_q(d,p)=\sum [0,1]\{v\in\ZZ^d\,:\,\|v\|_q\leq p\,,\ \gcd(v)=1\,,\ v\succ 0\}.$$
In other words, $H_q(d,p)$ is the Minkowski sum of the generators of $Z_q(d,p)$. We also consider the {\em positive primitive  zonotope} $Z^+_q(d,p)$ defined as the zonotope generated by the primitive integer vectors of $q$-norm at most $p$ with nonnegative coordinates:
$$Z^+_q(d,p)=\sum [-1,1]\{v\in\ZZ^d_+\,:\,\|v\|_q\leq p\,,\ \gcd(v)=1\}$$
where $\ZZ_+=\{0,1,\dots\}$. Similarly, we consider the Minkowski sum of the generators of $Z^+_q(d,p)$:
$$H^+_q(d,p)=\sum [0,1]\{v\in\ZZ^d_+\,:\,\|v\|_q\leq p\,,\ \gcd(v)=1\}.$$

\noindent
We illustrate the primitive  zonotopes with a few examples:
\begin{itemize}
\item[$(i)$]
For finite $q$, $Z_q(d,1)$ is generated by the $d$ unit vectors and forms  the  $\{-1,1\}^d$-cube, and  $H_q(d,1)$ is the $\{0,1\}^d$-cube.
\item[$(ii)$]
 $Z_1(d,2)$ is the permutahedron of type $B_d$ and thus, $H_1(d,2)$ is, up to translation, a lattice $(d,2d-1)$-polytope with $2^dd!$ vertices and diameter $d^2$. For example,
 $Z_1(2,2)$ is generated by $\{(0,1),(1,0),(1,1),(1,-1)\}$ and forms the octagon whose vertices are $\{(-3,-1),(-3,1),(-1,3),(1,3),(3,1),(3,-1),(1,-3),(-1,-3)\}$. $H_1(2,2)$ is, up to translation, a lattice $(2,3)$-polygon. 
$Z_1(3,2)$ is 
congruent to the  truncated cuboctahedron -- which is also called  great rhombicuboctahedron  -- 
and is the Minkowski sum of an octahedron and a cuboctahedron, see for instance Eppstein~\cite{zo}.  $H_1(3,2)$ is, up to translation, a lattice $(3,5)$-polytope with diameter $9$ and $48$ vertices.
\item[$(iii)$]
 $H^+_1(d,2)$ is the Minkowski sum of the permutahedron with the $\{0,1\}^d$-cube. Thus, $H^+_1(d,2)$ is a lattice $(d,d)$-polytope with  diameter ${d+1 \choose 2}$.
 \item[$(iv)$]
$Z_\infty(3,1)$  is congruent to the truncated small rhombicuboctahedron, see Figure~\ref{Einfty31} for an illustration, 
which is the Minkowski sum of a cube, a  truncated octahedron, and a rhombic dodecahedron, see for instance Eppstein~\cite{zo}. $H_\infty(3,1)$ is, up to translation, a lattice $(3,9)$-polytope with diameter $13$ and $96$ vertices.
 \begin{figure}[htb]
 \begin{center}
 \includegraphics[height=0.24\textheight]{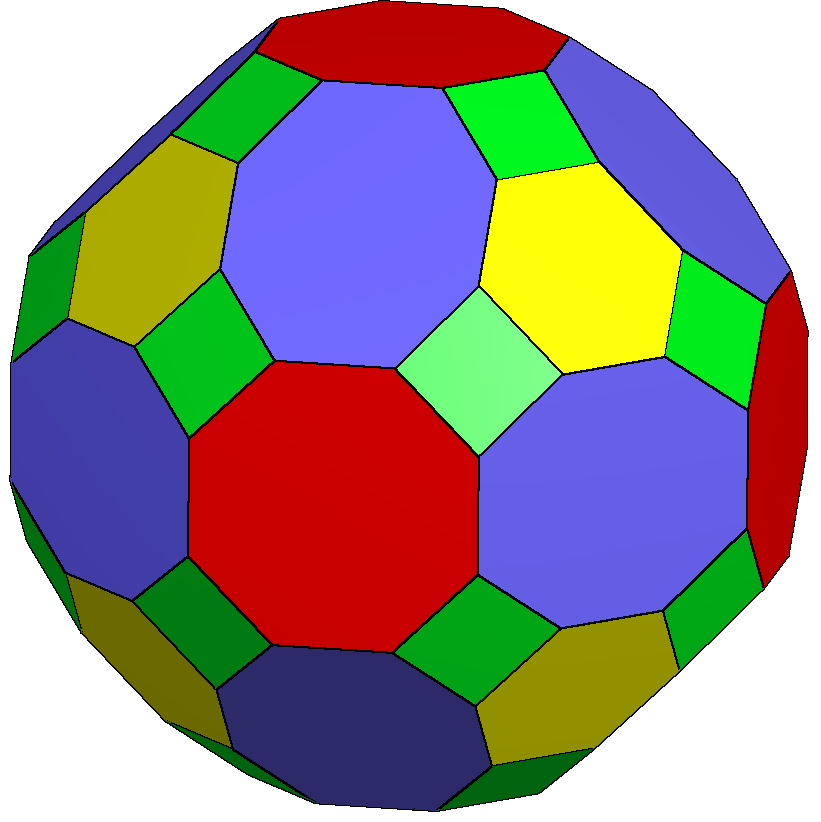}
   \caption{$Z_\infty(3,1)$ is congruent  to the truncated small rhombicuboctahedron}\label{Einfty31}
  \end{center}
 \end{figure}
\item[$(v)$]
$Z^+_\infty(2,2)$  is generated by $\{(0,1),(1,0),(1,1),(1,2),(2,1)\}$ and forms the decagon whose vertices are $\{(-5,-5),(-5,-3),(-3,-5),(-3,1),(-1,3),(1,-3),(3,-1),(3,5),$ $(5,3),$ $(5,5)\}$. $H^+_\infty(2,2)$
is a lattice $(2,5)$-polygon. 
 \end{itemize}

\subsection{Combinatorial properties}
We provide properties concerning $Z_q(d,p)$  and $Z^+_q(d,p)$, and in particular their  symmetry group, diameter, and vertices.  $Z_1(d,2)$ is the permutahedron of type $B_d$ as its generators form the root system of type $B_d$, see~\cite{H90}. Thus, $Z_1(d,2)$  has $2^dd!$ vertices  and  its symmetry group is $B_d$.  The properties listed in this section are extensions to $Z_q(d,p)$ of  known properties of $Z_1(d,2)$, and thus given without proof. We refer to Fukuda~\cite{F15}, Gr\"unbaum~\cite{G03}, and Ziegler~\cite{Z95} for polytopes and, in particular, zonotopes.
\begin{property}\label{sym}$\:$
\begin{itemize}
\item[$(i)$]
 $Z_q(d,p)$ is invariant under the symmetries induced by coordinate permutations and the reflections induced by sign flips.
\item[$(ii)$]
The sum $\sigma_q(d,p)$ of all the generators of $Z_q(d,p)$ is a vertex of both $Z_q(d,p)$ and $H_q(d,p)$. The origin is a vertex of $H_q(d,p)$ and $-\sigma_q(d,p)$ is a vertex of $Z_q(d,p)$.
\item[$(iii)$]
The coordinates of the vertices of  $Z_q(d,p)$ are odd, and thus the number of vertices of  $Z_q(d,p)$ is a multiple of $2^d$.
\item[$(iv)$]
$H_q(d,p)$ is, up to translation,  a lattice $(d,k)$-polytope where $k$ is the sum of the first coordinates of all generators of $Z_q(d,p)$.
\item[$(v)$]
The diameter of $Z_q(d,p)$, respectively $Z^+_q(d,p)$, is equal to the number of its generators.
\end{itemize}
\end{property}
\begin{property}\label{sym+}$\:$
\begin{itemize}
\item[$(i)$]
$Z^+_q(d,p)$ is centrally symmetric and invariant under the symmetries induced by coordinate permutations.
\item[$(ii)$]
The sum $\sigma^+_q(d,p)$ of all the generators of $Z^+_q(d,p)$ is a vertex of both $Z^+_q(d,p)$ and $H^+_q(d,p)$. The origin is a vertex of $H^+_q(d,p)$ and $-\sigma^+_q(d,p)$ is a vertex of $Z^+_q(d,p)$.
\end{itemize}
\end{property}
\noindent
A vertex $v$ of  $Z_q(d,p)$ is called {\em canonical} if $v_1\geq\dots\geq v_d>0$. Property~\ref{sym} item $(i)$ implies that the vertices of $Z_q(d,p)$ are all the coordinate permutations and sign flips of its canonical vertices.
\begin{property}\label{canonical}$\:$
\begin{itemize}
\item[$(i)$]
A canonical vertex $v$ of  $Z_q(d,p)$ is the unique maximizer of $\{\max c^Tx : x\in Z_q(d,p)\}$ for some vector $c$ satisfying $c_1>c_2>\dots>c_d>0$.
\item[$(ii)$]
$Z_1(d,2)$ has $2^dd!$ vertices corresponding to all coordinate permutations and sign flips of  the unique canonical vertex $\sigma_1(d,2)=(2d-1,2d-3,\dots,1)$.
\item[$(iii)$]
$Z^+_\infty(d,1)$ has at least  $2+2d!$ vertices including the $2d!$ permutations of $\pm\sigma(d)$ where $\sigma(d)$ is a vertex with pairwise distinct coordinates, and the 2 vertices $\pm\sigma^+_\infty(d,1)$.
\end{itemize}
\end{property}

\noindent
Enumerative questions concerning $H_q(d,p)$ and $H_q^+(d,p)$ have been studied in various settings. For example, the number of vertices of $H_\infty^+(d,1)$ corresponds to the OEI sequence A034997 giving the number of generalized retarded functions in quantum field theory, and the number of vertices of $H_\infty(d,1)$, which is the number of regions of hyperplane arrangements with $\{-1,0.1\}$-valued normals in dimension $d$, corresponds, up to a factor of  $2^dd!$,  to the OEI sequence A009997, see~\cite{OEI}  and references therein.

\section{Large diameter}\label{deltaZ}
 Let $\delta(d,k)$ be the maximum possible edge-diameter over all lattice $(d,k)$-polytopes. Naddef~\cite{N89} showed in 1989 that $\delta(d,1)=d$, Kleinschmidt and Onn~\cite{KO92} generalized this result in 1992 showing that $\delta(d,k)\leq kd$. In 2016, Del Pia and Michini~\cite{DM16} strengthened the upper bound to $\delta(d,k) \leq kd - \lceil d/2\rceil$ for $k\geq 2$, and showed that  $\delta(d,2)=\lfloor3d/2\rfloor$. Pursuing Del Pia and Michini's approach, Deza and Pournin~\cite{DP16} showed that $\delta(d,k) \leq kd - \lceil 2d/3\rceil-(k-3)$ for $k\geq 3$, and that $\delta(4,3)=8$.  Del Pia and Michini conclude their paper noting that the current lower bound for $\delta(d,k)$ is of order $k^{2/3}d$ and ask whether the gap between the lower and upper bounds could be closed, or at least reduced. The order $k^{2/3}d$ lower bound for $\delta(d,k)$ is a direct consequence of the determination of $\delta(2,k)$ which was investigated independently in the early nineties by Thiele~\cite{T91}, Balog and B\'ar\'any~\cite{BB91}, and Acketa and \v{Z}uni\'{c}~\cite{AZ95}. In this section, we highlight that $H_1(2,p)$ is the unique polygon achieving  $\delta(2,k)$ for a proper $k$, and that a Minkowski sum of a proper subset of the generators of  $H_1(d,2)$ achieves a diameter of $\lfloor (k+1)d/2\rfloor$ for all $k\leq 2d-1$.
 \subsection{$H_1(2,p)$  as  a lattice polygon with large diameter}

Finding lattice polygons with the largest diameter; that is, to determine $\delta(2,k)$, was investigated independently in the early nineties by Thiele~\cite{T91}, Balog and B\'ar\'any~\cite{BB91}, and Acketa and \v{Z}uni\'{c}~\cite{AZ95}. This question can be found in Ziegler's book~\cite{Z95} as Exercise 4.15.  The answer is summarized in Proposition~\ref{Z12p} where $\phi(j)$ is the Euler totient function counting positive integers less than or equal to $j$ and relatively prime with $j$. Note that $\phi(1)$ is set to 1.
\begin{proposition}\label{Z12p}
$H_1(2,p)$ is, up to translation, a lattice $(2,k)$-polygon with $k=\sum_{j=1}^{p} j\phi(j)$ where $\phi(j)$ denotes the Euler totient function. The diameter of $H_1(2,p)$ is $2\sum_{j=1}^{p} \phi(j)$ and satisfies $\delta(H_1(2,p))=\delta(2,k)$. Thus, $\delta(2,k)=6(\frac{k}{2\pi})^{2/3}+O(k^{1/3}\log k)$.
\end{proposition}
\noindent
Note that lattice polygons can be associated to set of integer-valued vectors adding to zero and such that no pair of vectors are positive multiples of each other. Such set of vectors forms a $(2,k)$-polygon with $2k$ being the maximum between the sum of the norms of the first coordinates of the vectors and the sum of the norms of the second coordinates of the vectors. Then, for $k=\sum_{j=1}^{p}  j\phi(j)$ for some $p$, one can show that $\delta(2,k)$ is achieved uniquely by a translation of $H_1(2,p)$. For $k\neq\sum_{j=1}^{p}  j\phi(j)$ for any $p$, $\delta(2,k)$ is achieved by a translation of a Minkowski sum of an appropriate subset of the generators of $H_1(2,p)$ including all generators of $H_1(2,p-1)$ for an appropriate $p$. For the order of  $\sum_{j=1}^{p} \phi(j)$, respectively $\sum_{j=1}^{p} j\phi(j)$, being $\frac{3p^2}{\pi^2} + O(p\ln p)$, respectively $\frac{2p^3}{\pi^2} + O(p^2\ln p)$, we refer to~\cite{3}. The first values of $\delta(2,k)$ are given in Table~\ref{TE}.
\begin{table}[htb]
\begin{center}
\begin{tabular}{c||c|c|c|c|c|c|c|c|c|c|c|c|c|c|c|c|c}
$p$ of $H_1(2,p)$ &1 &  & 2 &  &  &  &  &  & 3 &  &   &  &  &  &  &  & 4  \\
 \hline
$k$ &1 & 2 & 3 & 4 & 5 & 6 & 7 & 8 & 9 & 10 & 11 & 12 & 13 & 14 & 15 & 16 & 17 \\
 \hline
$\delta(2,k)$ &2 & 3&4 &4& 5& 6& 6&7 &8& 8& 8 & 9&10 & 10&10 &11 &12
\end{tabular}
\caption{Relation between $H_1(2,p)$ and $\delta(2,k)$}\label{TE}
   \end{center}
\end{table}

 \subsection{$H_1(d,2)$  as a lattice polytope with large diameter}
As pointed out by Vincent Pilaud, a lower bound of $kd/2$ for $\delta(d,k)$ for appropriate $k<d$ can be achieved by considering a graphical zonotope $H_{\cal G}$; that is, the Minkowski sum of the line segments $[e_i,e_j]$ for all edges $ij$ of a given graph ${\cal G}$.  Consider the graphical zonotope $H_{{\cal C}(d,k)}$ associated to the circulant graph ${\cal C}(d,k)$ of degree $k$ on $d$ nodes.  One can check that $H_{{\cal C}(d,k)}$ is a lattice $(d,k)$-polytope with diameter $kd/2$. In this section, pursuing this approach, we  show that a Minkowski sum of a proper subset of the generators of  $H_1(d,2)$ yields $\delta(d,k)\geq \lfloor (k+1)d/2\rfloor$  for all  $k\leq 2d-1$.
\begin{theorem}\label{Z1d2}
For $k\leq 2d-1$, there is a subset  of the generators of  $H_1(d,2)$ whose Minkowski sum is, up to translation, a lattice $(d,k)$-polytope with diameter $\lfloor (k+1)d/2\rfloor$. So for  $k\leq 2d-1$ we have $\delta(d,k)\geq \lfloor (k+1)d/2\rfloor$. For instance, 
$H_1^+(d,2)$ is a lattice $(d,d)$-polytope with diameter $d+1\choose 2$, and
$H_1(d,2)$ is, up to translation, a lattice $(d,2d-1)$-polytope with diameter $d^2$.
\begin{proof}
We first note that the number of generators of $H_1(d,2)$ is $d^2$. The generators of $H_1(d,2)$ are $\{-1,0,1\}$-valued $d$-tuples: $d$ permutations of $(1,0,\dots,0)$, ${d \choose 2}$ permutations of $(1,1,0,\dots,0)$, and ${d \choose 2}$ permutations of $(1,-1,0,\dots,0)$. Thus, $\delta(H_1(d,2))=d^2$ by Property~\ref{sym} item $(v)$. As the sum of the first coordinates of the generators of $H_1(d,2)$ is $2d-1$, $H_1(d,2)$ is, up to translation, a lattice $(d,2d-1)$-polytope by  Property~\ref{sym} item $(iv)$.
Consider first the case when $d$ is even.  The first $d-1$ subsets are obtained by removing  from the current subset of generators  of  $H_1(d,2)$  a set of $d/2$ generators taken among the ${d \choose 2}$ permutations of $(1,-1,0,\dots,0)$. The removed $d-1$ subsets correspond to $d-1$ disjoint perfect matchings of the complete graph $K_d$ where the nonzero $i^{th}$ and $j^{th}$  coordinates of a generator $(\dots,1,\dots,-1,\dots)$ correspond to the edge $[i,j]$. The first  perfect matching is $[1,2],[3,d],[4,d-1],\dots,[d/2+1,d/2+2]$. The next perfect matching is obtained by changing $d$ to 2, and $i$ to $i+1$ for all other entries except $1$, which remains unchanged. This procedure yields $d-1$ disjoint perfect matchings as, placing the vertices $2$ to $d$ on a circle around $1$ where the edge $[1,2]$ is vertical and the edges $[3,d],[4,d-1],\dots,[d/2+1,d/2+2]$ are horizontal, the procedure corresponds to the $d-1$ rotations of the initial perfect matching, see~\cite[Chapter~12]{B83}. As these $d-1$ perfect matchings correspond to all the generators of $H_1(d,2)$ which are permutations of $(1,-1,0,\dots,0)$, the procedure ends with a subset of the generators  of  $H_1(d,2)$ forming the $d+1\choose 2$ generators of  $H_1^+(d,2)$. We can then repeat the same procedure where the nonzero $i^{th}$ and $j^{th}$  coordinates of a generator $(\dots,1,\dots,1,\dots)$ correspond to the edge $[i,j]$ of $K_d$, and similarly obtain $d-1$ disjoint perfect matchings. The procedure now ends  with a subset of the generators  of  $H_1(d,2)$  forming $H_1(d,1)$; that is the unit cube. One can check that if the Minkowski sum $H$ of the current subset of generators  of  $H_1(d,2)$  is a lattice $(d,k)$-polytope of diameter $\delta(H)$, removing the $d/2$ generators corresponding to a perfect matching yields a lattice $(d,k-1)$-polytope of diameter $\delta(H)-d/2$. Thus, starting from $H_1(d,2)$  which is a  $(d,2d-1)$-polytope with diameter $d^2$, we obtain a  $(d,k)$-polytope with diameter $(k+1)d/2$ for all $k\leq 2d-1$.
The case when $d$ is odd is similar.  The removed subsets are of alternating sizes  $\lceil d/2\rceil$ and  $\lfloor d/2\rfloor$. Adding a dummy vertex $d+1$ to $K_d$, we consider the $d$ disjoint perfect matching of $K_{d+1}$ described for even $d$. The first subset consists of the $\lceil d/2\rceil$ edges where $[3,d+1]$ is replaced by $[3,5]$, the second subset  consists of the $\lfloor d/2\rfloor$ edges where $[5,d+1]$ is removed, the third subset consists of the $\lceil d/2\rceil$ edges where $[7,d+1]$ is replaced by $[7,9]$, and so forth. As for even $d$, one can check that if the Minkowski sum $H$ of the current subset of generators  of  $H_1(d,2)$  is a lattice $(d,k)$-polytope of diameter $\delta(H)$, removing the described $\lceil d/2\rceil$, respectively  $\lfloor d/2\rfloor$,  generators yields a lattice $(d,k-1)$-polytope of diameter $\delta(H)-\lceil d/2\rceil$,  respectively $\delta(H)-\lfloor d/2\rfloor$. Thus, starting from $H_1(d,2)$  which is a  $(d,2d-1)$-polytope with diameter $d^2$, we obtain a  $(d,k)$-polytope with diameter $\lfloor (k+1)d/2\rfloor$ for all $k\leq 2d-1$.
\end{proof}
\end{theorem}
\begin{conjecture}\label{C3}
$\delta(d,k)\leq \lfloor (k+1)d/2 \rfloor$, and $\delta(d,k)$ is achieved, up to translation, by a Minkowski sum of lattice vectors.
\end{conjecture}
 \noindent
Note that  Conjecture~\ref{C3} holds for all known values of $\delta(d,k)$ given in Table~\ref{delta(d.k)}, and hypothesizes, in particular,  that $\delta(d,3)=2d$.
Note that $\delta(d,3)=2d$ for $d\leq 4$, $2d\leq \delta(d,3)\leq\lfloor7d/3\rfloor-1$ when $d\not\equiv2\mod{3}$, and $\delta(d,3)\leq\lfloor7d/3\rfloor$ when $d\equiv2\mod{3}$, see~\cite{DP16}.

\begin{table}[bht]
\begin{center}
  \begin{tabular}{cc|ccccccccccc}
   & & \multicolumn{9}{|c}{$k$} \\
   & $\delta(d,k)$ & 1 & 2 & 3 & 4 & 5 & 6 & 7 & 8 & 9 & 10\\
   \hline
   \multirow{6}{*}{$d$} & 1 & 1 &  1 & 1  & 1 & 1 & 1 & 1 & 1 & 1 & 1\\
   & 2 & 2 &  3 & 4  & 4 & 5  & 6 & 6 & 7 & 8 & 8 \\
   & 3 & 3 &4 &  6 & \\
   & 4 & 4 & 6 &  8 & \\
   & $\vdots$ & $\vdots$ & $\vdots$ &  & \\
   & $d$ & $d$ & $\lfloor 3d/2\rfloor$ &  & \\
  \end{tabular} \caption{Largest diameter $\delta(d,k)$  over all lattice $(d,k)$-polytopes}\label{delta(d.k)}
\end{center}
\end{table}

\noindent
Soprunov and Soprunova~\cite{SS16} considered  the Minkowski length $L(P)$ of a lattice polytope $P$; that is, the largest number of lattice segments whose Minkowski sum is contained in $P$. Considering the special case when $P$ is the $\{0,k\}^d$-cube, let $L(d,k)$ denote the Minkowski length of $\{0,k\}^d$-cube. For example, the Minkowski length of the $\{0,1\}^d$-cube satisfies  $L(d,1)=d$. One can check that the generators of $H_1(d,2)$ form the largest, and unique, set of primitive lattice vectors which Minkowski sum fits within the $\{0,k\}^d$-cube for $k=2d-1$; that is, for $k$ being the sum of the first coordinates of the $d^2$ generators of $H_1(d,2)$. Thus, $L(d,2d-1)=\delta(H_1(d,2))=d^2$. Similarly, the constructions used in Proposition~\ref{Z12p} and~\ref{Z1d2} imply that $L(2,k)=\delta(2,k)$ and $L(d,k)=\lfloor (k+1)d/2\rfloor$ for $k\leq 2d-1$.

\section{Multicriteria matroid optimization}\label{vZ}
We consider the convex multicriteria matroid optimization framework of Melamed, Onn and Rothblum in~\cite{MO,Onn,OR}, and show that $H_\infty(d,p)$ settles its computational complexity.

Call $S\subset\{0,1\}^n$ a {\em matroid} if it is the set of the indicators of bases of a matroid over $\{1,\dots,n\}$. For instance, $S$ can be the set of indicators of spanning trees
in a connected graph with $n$ edges. We assume that the matroid is presented by an
{\em independence oracle} that, queried on $y\in\{0,1\}^n$, asserts whether or not
$y\leq x$ for some $x\in S$. The standard linear optimization problem over the matroid $S$ is: given a {\em utility vector} $w\in\ZZ^n$, find a basis which maximizes the utility $wx$,
that is, solve $\{\max wx:x\in S\}$. This problem is well known to be easily solvable by the
{\em greedy algorithm}. Generalizing this problem to $d$ criteria, we are given a $d\times n$ integer {\em utility matrix} $W$ whose $i^{th}$ row $W_i$ gives the utility $W_ix$ of basis $x\in S$ under criterion $i$, so the vector $Wx\in\ZZ^d$ represents the $d$ utility values of basis $x$ under the $d$ criteria. These values are then traded-in by a convex function $f:\RR^d\rightarrow\RR$. We assume that $f$ is presented by a {\em comparison oracle} that, queried on vectors $x,y\in\ZZ^d$, asserts whether or not $f(x)<f(y)$.
The {\em multicriteria matroid optimization problem} is then: find a basis which maximizes
the traded-in utility $f(Wx)$; that is, solve $\{\max f(Wx):x\in S\}$,
making use of the oracle presentations of $S$ and $f$.

Let $\conv(WS)=\conv\{Wx:x\in S\}$ be the projection to $\RR^d$ of $\conv(S)$ by $W$.
As detailed in~\cite[Chapter 2]{Onn}, the projection polytope $\conv(WS)$ and its vertices play a key role in solving our optimization problem, since for any convex function $f$ there is an optimal solution $x\in S$ whose projection $u=Wx$ is a vertex of $\conv(WS)$. Thus, the
convex multicriteria problem can be solved by enumerating the set of vertices of $\conv(WS)$, picking a vertex $u$ attaining a maximum value $f(u)$, and finding $x\in S$ with $Wx=u$.
However, direct computation of $\conv(WS)$ and enumeration of its vertices
are intractable since typically $S$ has exponentially many points.

Following~\cite{MO}, we consider nonnegative utilities, so that, for some positive integer $p$,  for all $i,j$, the utility $W_{i,j}$ of element $j$ of the ground set of the matroid under criterion $i$ is in $\{0,1,\dots,p\}$. We call such utility matrices
{\em $p$-bounded}. Let $m(d,p)$ be the number of vertices of $H_\infty(d,p)$.
Theorem~\ref{Matroid_Theorem} settles the computational complexity of the multicriteria optimization problem by showing that the maximum number of vertices of the projection polytope $\conv(WS)$ of any matroid $S$ on $n$ elements and any $d$-criteria $p$-bounded utility matrix; that is, $W\in\{0,1,\dots,p\}^{d\times n}$, is equal to $m(d,p)$, and hence is in particular
{\em independent} of $n,S$, and $W$. Below we use the following. The {\em normal cone} of polytope $P\subset\RR^n$ at its vertex $v$ is the relatively open cone of vectors $h\in\RR^n$ such that $v$ is the unique maximizer of $hx$ over $P$. A polytope $H$ {\em refines} a polytope $P$ if the normal cone of $H$ at every vertex of $H$ is contained in the normal cone of $P$ at some vertex of $P$. Then, the closure of each normal cone of $P$ is the union of closures
of normal cones of $H$ and $P$ has no more vertices than $H$.

\begin{theorem}\label{Matroid_Theorem}
Let $d,p$ be any positive integers. Then, for any positive integer $n$, any matroid $S\subset\{0,1\}^n$, and any $d$-criteria $p$-bounded utility matrix $W$, the primitive
zonotope $H_\infty(d,p)$ refines $\conv(WS)$. Moreover, $H_\infty(d,p)$ is a translation of
$\conv(WS)$ for some matroid $S$ and $d$-criteria $p$-bounded utility matrix $W$. Thus, the maximum number of vertices of $\conv(WS)$ for any $n$, any matroid $S\subset\{0,1\}^n$, and any $d$-criteria $p$-bounded utility matrix $W$, equals the number $m(d,p)$ of vertices of $H_\infty(d,p)$, and hence is in particular independent of $n,S$, and $W$.
Also, for any fixed $d$ and convex $f:\RR^d\rightarrow\RR$, the multicriteria matroid optimization problem can be solved using a number of arithmetic operations and queries
to the oracles of $S$ and $f$ which is polynomial in $n$ and $p$ using
$m(d,p)$ greedily solvable linear matroid optimization counterparts.

\begin{proof}
First we show that $H_\infty(d,p)$ refines $\conv(WS)$ for every matroid $S$. It is known that if $G$ is a finite set of vectors such that every edge in a polytope $P$ is parallel to some $g\in G$ then the zonotope $H=\sum[0,1]G$ refines $P$, see~\cite{GS,OR}. Now, for any matroid $S\subset\{0,1\}^n$, any edge of $\conv(S)$ is parallel to the difference ${\bf 1}_i-{\bf 1}_j$ between a pair of unit vectors in $\RR^n$, see~\cite[Chapter 2]{Onn}. Therefore any edge of the projection $\conv(WS)$ is parallel to the difference $W^i-W^j$ between a pair of columns of $W$. Since $W^i,W^j\in\{0,1,\dots, p\}^d$, we have that $W^i-W^j\in\{0,\pm1,\dots,\pm p\}^d$,
and so it follows that every edge of $\conv(WS)$ is parallel to some vector in
$$G(d,p)\ =\ \{v\in\ZZ^d\,:\,\|v\|_\infty\leq p\,,\ \gcd(v)=1\,,\ v\succ 0\}\ .$$
It now follows that the primitive zonotope $H_\infty(d,p)=\sum[0,1]G(d,P)$ refines $\conv(WS)$.
Next we construct a matroid $S$ and $d$-criteria $p$-bounded utility matrix such that $H_\infty(d,p)$ is a translation of $\conv(WS)$.
For vector $g\in\ZZ^d$ let $g^+,g^-\in\ZZ^d_+$ be its
{\em positive and negative parts}; that is, the unique nonnegative vectors with disjoint support satisfying $g=g^+-g^-$.
Let $r=|G(d,p)|$ and $n=2r$. Order the vectors $G(d,p)$ arbitrarily $G(d,p)=\{g_1,\dots,g_r\}$. Consider the graph which is a path of length $r$
whose edges are labeled by $g_1,\dots,g_r$. Now replace each edge $g_i$ by two
parallel edges labeled by $g^+_i$ and $g^-_i$. Let $S\subset\{0,1\}^n$ be the graphic matroid of the resulting graph on $n$ edges. Let $W$ be the $d$-criteria $p$-bounded matrix
$W=[g^+_1,g^-_1,\dots,g^+_r,g^-_r]$.
Now, there is a bijection between bases $x\in S$ and subsets $I\subseteq\{1,\dots,r\}$  where we put $i$ in $I$ if $x$ chooses $g^+_i$ from the parallel pair $\{g^+_i,g^-_i\}$.
Then, for corresponding $x$ and $I$:
$$Wx\ =\ \sum_{i\in I} g^+_i+\sum_{i\notin I} g^-_i\ = \
\sum_{i=1}^r g^-_i+\sum_{i\in I}(g^+_i-g^-_i)+\sum_{i\notin I}(g^-_i-g^-_i)\ =\
\sum_{i=1}^r g^-_i+\sum_{i\in I}g\ .$$
Thus,
\begin{eqnarray*}
\conv(Wx\,:\,x\in S\}& = &
\conv\left(\sum_{i=1}^r g^-_i+\sum_{i\in I}g\,:\,I\subseteq\{1,\dots,r\}\right)\\
& = & \sum_{i=1}^r g^-_i+\conv\left(\sum G\,:\,G\subseteq G(d,p)\right)
\ =\ \sum_{i=1}^r g^-_i+H_\infty(d,p)\ .
\end{eqnarray*}

\noindent
We now conclude that this implies that the maximum number of vertices of $\conv(WS)$ for any $n$, any matroid $S\subset\{0,1\}^n$, and any $d$-criteria $p$-bounded utility matrix $W$,
equals the number $m(d,p)$ of vertices of $H_\infty(d,p)$. Indeed, for every $W$ and $S$, the number of vertices of $\conv(WS)$ is at most $m(d,p)$ since $H_\infty(d,p)$ refines $\conv(WS)$. In addition,  $m(d,p)$ is attained as the number of vertices of the translation $\conv(WS)$ of $H_\infty(d,p)$ for  $S$ and $W$ constructed above.

We proceed with the statements about the algorithmic complexity of the multicriteria matroid optimization problem. Let $d$ be fixed, and  $g(d,p)=|G(d,p)|=O(p^d)$ be the number of generators of $H_\infty(d,p)$.
Then, as detailed in~\cite[Chapter 2]{Onn}, the number of vertices of $H_\infty(d,p)$ satisfies $m(d,p)=O(g(d,p)^{d})=O(p^{d^2})$
which is polynomial in $p$, and in that much time, all vertices of $H_\infty(d,p)$
can be enumerated along with, for each vertex $v\in\ZZ^d$,  a vector $h_v$ in the normal cone of $H_\infty(d,p)$ at $v$. This preprocessing depends only on $p$.

\ignore{
 \textcolor{red}{
We proceed with the statements about the algorithmic complexity of the multicriteria matroid optimization problem. Let $d$ be fixed, and  $g(d,p)=|G(d,p)|=O(p^d)$ be the number of generators of $H_\infty(d,p)$.
Then, as detailed in~\cite[Chapter 2]{Onn}), the number of vertices of $H_\infty(d,p)$ satisfies $m(d,p)=O(g(d,p)^{d})=O(p^{d^2})$; and thus polynomial in $p$.
Consequently, all the vertices $v$ of $H_\infty(d,p)$, along with a vector $h_v$ in the normal cone of $H_\infty(d,p)$ at $v$, can be enumerated in polynomial time.
}
ignore}

Now, let $n$, matroid $S\subset\{0,1\}^n$, and $d$-criteria $p$-bounded utility matrix $W$ be given. For each of the $m(d,p)=O(p^{d^2})$ vertices $v$ of $H_\infty(d,p)$ we solve the standard linear optimization problem over $S$ with utility vector $w_v=h_vW\in\ZZ^n$ by the greedy algorithm using the independence oracle of $S$ and find an optimal basis $x_v\in S$.
We collect the projections $Wx_v$ of all these optimal bases $x_v$ corresponding to the vertices
$v$ of $H_\infty(d,p)$ in a set $U\subset\ZZ^d$. We now claim that every vertex $u$ of $\conv(WS)$ lies in $U$. Consider such a vertex $u$ and let $x\in S$ be such that $u=Wx$.
Since $H_\infty(d,p)$ refines $\conv(WS)$, there is a vertex $v$ of $H$ such that the normal cone of $H$ at $v$ is contained in the normal cone of $\conv(WS)$ at $u$. Therefore, $h_v$ is maximized over $\conv(WS)$ uniquely at $u=Wx$. We claim that $u=Wx_v\in U$. Assume that not, then we get $w_vx=h_vWx>h_vWx_v=w_vx_v$. hence a contradiction. 
Thus, we find a vertex $v$ of $H_\infty(d,p)$ such that $u=Wx_v$ maximizes $f(u)$ over $U$ using the comparison oracle of $f$, and conclude that $x_v\in S$ is the optimal solution to the multicriteria matroid problem.
\end{proof}
\end{theorem}

\begin{example}\label{multicriteria_example}
{\rm Let $(d,p)=(2,1)$, we describe $H_\infty(2,1)$, the matroid $S$, and  the matrix $W$
constructed in the proof of Theorem~\ref{Matroid_Theorem} such that $H_\infty(2,1)$ is a translation of $\conv(WS)$:
$G(2,1)=\{(0,1),(1,-1),(1,0),(1,1)\}$, $r=4$, $n=8$, and 
$$S\ =\ \{x\in\{0,1\}^8\,:\,x_{2i-1}+x_{2i}=1\,,\ \ i=1,2,3,4\}\ ,\quad
W\ =\
\left(
\begin{array}{cccccccc}
  0 & 0 & 1 & 0 & 1 & 0 & 1 & 0 \\
  1 & 0 & 0 & 1 & 0 & 0 & 1 & 0 \\
\end{array}
\right)\ .$$
See Figure~\ref{octagon} for an illustration of $\conv(WS)=(0,1)+H_\infty(2,1)$ with its $m(2,1)=8$ vertices, and a vector $h_v$
in the normal cone of $H_\infty(2,1)$ at each vertex $v$ of $\conv(WS)$.\\

\noindent
Now consider the following bicriteria matroid optimization problem with
$f=\|\cdot\|_2^2$; that is, $f(y_1,y_2)=y_1^2+y_2^2$, over a uniform matroid
and $2$-criteria $1$-bounded utility matrix given by
$$U_{12}^6=\{x\in\{0,1\}^{12}\,:\,\sum_{i=1}^{12}x_i=6\}\ ,\quad
W\ =\
\left(
\begin{array}{cccccccccccc}
  0 & 0 & 1 & 1 & 0 & 0 & 1 & 1 & 0 & 0 & 1 & 1\\
  0 & 1 & 0 & 1 & 0 & 1 & 0 & 1 & 0 & 1 & 0 & 1\\
\end{array}
\right)\ .$$
We solve the problem using the algorithm provided in the proof of Theorem~\ref{Matroid_Theorem}. For each of the $m(2,1)=8$ vectors $h_v$ in normal cones in Figure~\ref{octagon},
we solve the linear optimization counterpart over $U_{12}^6$ with utility vector $w_v=h_vW\in\ZZ^{12}$ by the greedy algorithm and find an optimal basis $x_v\in U_{12}^6$.
For instance, for $h_v=(1,2)$, greedily found  $w_v$ and $x_v$ are:
$$w_v =
\left(
\begin{array}{cccccccccccc}
  0 & 2 & 1 & 3 & 0 & 2 & 1 & 3 & 0 & 2 & 1 & 3
\end{array}
\right)\,,\quad
x_v =
\left(
\begin{array}{cccccccccccc}
  0 & 1 & 0 & 1 & 0 & 1 & 0 & 1 & 0 & 1 & 0 & 1
\end{array}
\right)\,.
$$
Then, $Wx_v=(3,6)$ so $x_v$ has objective value $f(Wx_v)=45$.
We repeat this for all $8$ vectors $h_v$ and find the best. Here, the $x_v$ above is indeed an optimal solution to the bicriteria problem.
Note that $U_{12}^6$ has ${12\choose 6}=924$ bases and for matroids on ground sets with larger $n$ the number of bases typically grows exponentially so solving our problem by exhaustive search is unreasonable. Instead, our algorithm solves any bicriteria $1$-bounded matroid problem in polynomial time by always greedily solving only $8$ linear counterparts,
each in time linear in $n$.
}
\begin{figure}[hbt]
\centerline{\includegraphics[scale=.4]{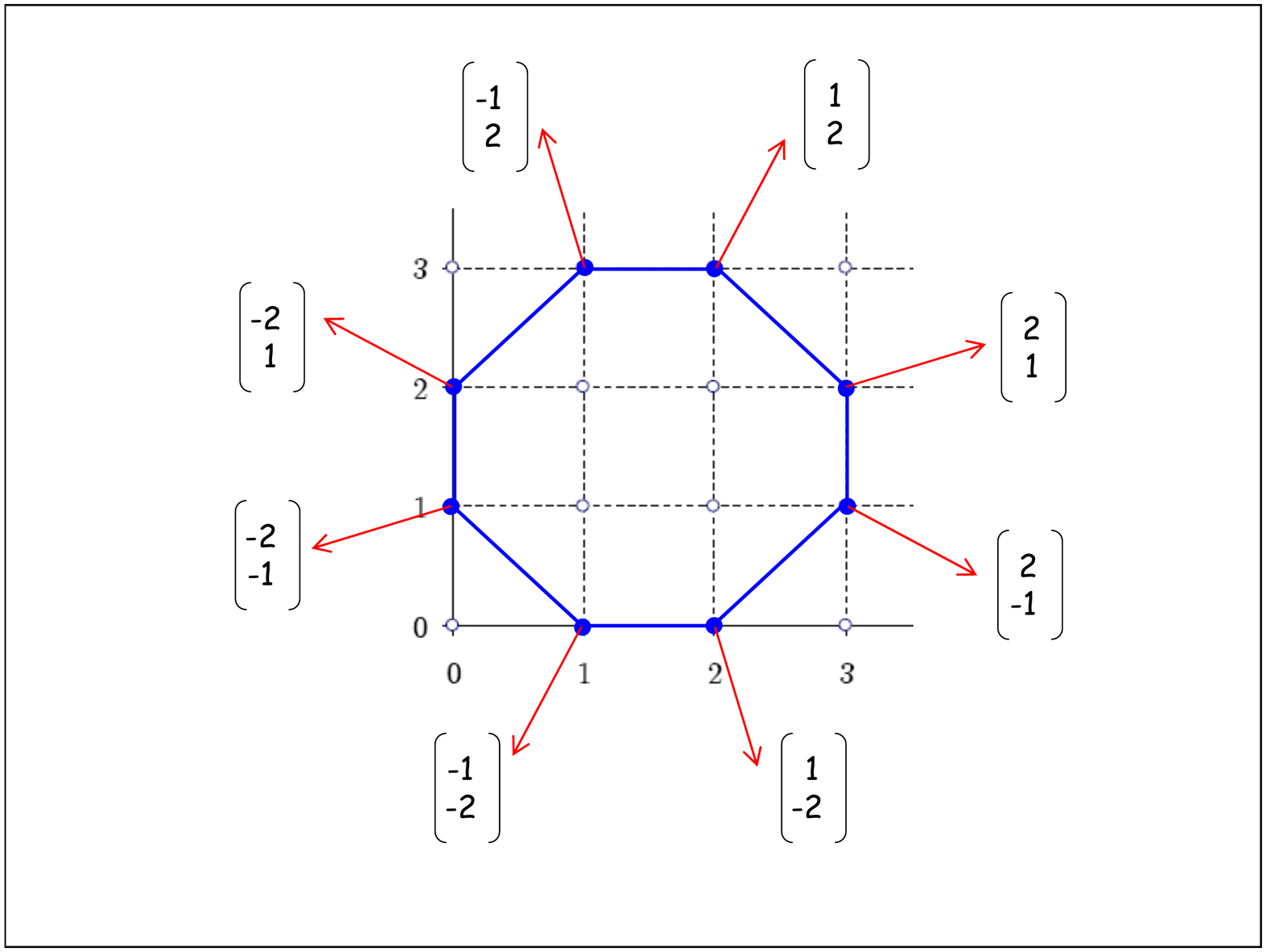}}
\caption{$\conv(WS)=(0,1)+H_\infty(2,1)$ and associated vectors $h_v$}
\label{octagon}
\end{figure}
\end{example}

\noindent
Theorem~\ref{Matroid_Theorem} asserts that the number $m(d,p)$ of vertices of $H_\infty(d,p)$ is the largest number of linear counterparts needed to solve any $d$-criteria $p$-bounded problem. While we do not know much about the exact value of $m(d,p)$, we examine $m(d,p)$ for small $d$ or $p$ in the rest of this section: $m(1,p)$ is trivially equal to $2$,  $m(2,p)$ is given in Proposition~\ref{dimension_2}, and the first values of $m(d,1)$ are given after Proposition~\ref{dimension_2}.

\begin{proposition}\label{dimension_2}
The number of vertices of $H_\infty(2,p)$ satisfies $m(2,p)=8\sum_{j=1}^p\phi(j)$ where $\phi(j)$ is the Euler totient function counting positive integers less than or equal to $j$ and relatively prime with $j$.
\begin{proof}
Since the generators of $H_q(d,p)$ are pairwise linearly independent, the diameter of $H_q(d,p)$ equals the number of generators. For $d=2$, the number of vertices of $H_q(2,p)$ is twice
the diameter. Thus, $m(2,p)$ is twice the number of generators of $H_{\infty}(2,p)$.
Now, $H_{\infty}(2,p)$ has $4\phi(1)$ generators $(1,0),(0,1),(1,1)$, $(1,-1)$. In addition, for $p\geq 2$,  $H_{\infty}(2,p)$ has the $4\phi(j)$ generators $(i,j)$, $(i,-j)$, $(j,i)$, $(j,-i)$ for $j=2,\dots,p$ where $i$ runs through all positive integers less than or equal to $j$ and relatively prime with $j$. So $H_{\infty}(2,p)$ has $4\sum_{j=1}^p\phi(j)$
generators and $m(2,p)=8\sum_{j=1}^p\phi(j)$ vertices.
\end{proof}
\end{proposition}

\noindent
For instance, the first values of $m(2,p)$ are:
$$m(2,1)=8,\ m(2,2)=16,\ m(2,3)=32,\ m(2,4)=48,\ m(2,5)=80,\ m(2,6)=96.$$
Turning to the number $m(d,1)$ of vertices of $H_{\infty}(d,1)$, no closed-form expression is known but, as mentioned at the end of Section~\ref{eulerZ}, $m(d,1)$ corresponds to the OEI sequence A009997, up to a factor of $2^dd!$, and the first values are:
$$m(2,1)=8,\ m(3,1)=96,\ m(4,1)=5\ 376,\ m(5,1)=1\ 981\ 440,\ m(6,1)=5\ 722\ 536\ 960.$$
To solve a $6$-criteria $1$-bounded matroid problem may require about $6$ billion  linear counterparts!

\section*{Acknowledgment}
The authors thank the anonymous referees, Johanne Cohen, Nathann Cohen, Komei Fukuda, and Aladin Virmaux for valuable comments and suggestions and for informing us of reference~\cite{T91},  Emo Welzl and G\"{u}nter Ziegler for helping us access Thorsten Thiele's Diplomarbeit, Dmitrii Pasechnik for pointing out reference~\cite{SS16} and the concept of Minkowski length, and Vincent Pilaud for pointing out graphical zonotopes and that $Z_1(d,2)$ is the permutahedron of type $B_d$. This work was partially supported by the Natural Sciences and Engineering Research Council of Canada Discovery Grant program (RGPIN-2015-06163), by the Digiteo Chair C\&O program, and by the Dresner Chair at the Technion.

\bibliographystyle{plain}
\bibliography{euler_arxiv}

\begin{thebibliography}{10}

\bibitem{AZ95}
Dragan Acketa and Jovi\v{s}a \v{Z}uni\'{c}.
\newblock On the maximal number of edges of convex digital polygons included
  into an $m\times{m}$-grid.
\newblock {\em Journal of Combinatorial Theory A}, 69:358--368, 1995.

\bibitem{ABGM14}
Xavier Allamigeon, Pascal Benchimol, St{\'e}phane Gaubert, and Michael Joswig.
\newblock Long and winding central paths.
\newblock {\em arXiv:1405.4161}, 2014.

\bibitem{BB91}
Antal Balog and Imre B\'{a}r\'{a}ny.
\newblock On the convex hull of the integer points in a disc.
\newblock In {\em Proceedings of the Seventh Annual Symposium on Computational
  Geometry}, pages 162--165, 1991.

\bibitem{B83}
Claude Berge.
\newblock {\em Graphes}.
\newblock Gauthier-Villars, 1983.

\bibitem{BSEHN14}
Nicolas Bonifas, Marco Di~Summa, Friedrich Eisenbrand, Nicolai H\"{a}hnle, and
  Martin Niemeier.
\newblock On sub-determinants and the diameter of polyhedra.
\newblock {\em Discrete and Computational Geometry}, 52:102--115, 2014.

\bibitem{BDF16}
Steffen {Borgwardt}, Jes{\'u}s {De Loera}, and Elisabeth {Finhold}.
\newblock The diameters of transportation polytopes satisfy the {H}irsch
  conjecture.
\newblock {\em arXiv:1603.00325}, 2016.

\bibitem{DM16}
Alberto {Del Pia} and Carla {Michini}.
\newblock On the diameter of lattice polytopes.
\newblock {\em Discrete and Computational Geometry}, 55:681--687, 2016.

\bibitem{DP16}
Antoine Deza and Lionel Pournin.
\newblock Improved bounds on the diameter of lattice polytopes.
\newblock {\em arXiv:1610.00341}, 2016.

\bibitem{zo}
David Eppstein.
\newblock Zonohedra and zonotopes.
\newblock {\em Mathematica in Education and Research}, 5:15--21, 1996.

\bibitem{F15}
Komei Fukuda.
\newblock Lecture notes: Polyhedral computation.
\newblock
  \url{http://www-oldurls.inf.ethz.ch/personal/fukudak/lect/pclect/notes2015/}.

\bibitem{GS}
Peter Gritzmann and Bernd Sturmfels.
\newblock Minkowski addition of polytopes: complexity and applications to
  {G}r\"obner bases.
\newblock {\em SIAM Journal on Discrete Mathematics}, 6:246--269, 1993.

\bibitem{G03}
Branko Gr{\"u}nbaum.
\newblock {\em Convex Polytopes}.
\newblock Graduate Texts in Mathematics. Springer, 2003.

\bibitem{3}
Godfrey Hardy, Edward Wright, David Hearth-Brown, and Joseph Silverman.
\newblock {\em An introduction to the theory of numbers}.
\newblock Clarendon Press Oxford, 1979.

\bibitem{H90}
James Humphreys.
\newblock {\em Reflection groups and Coxeter groups}.
\newblock Cambridge Studies in Advanced Mathematics. Cambridge University
  Press, 1990.

\bibitem{KK92}
Gil Kalai and Daniel Kleitman.
\newblock A quasi-polynomial bound for the diameter of graphs of polyhedra.
\newblock {\em Bulletin of the American Mathematical Society}, 26:315--316,
  1992.

\bibitem{KO92}
Peter Kleinschmidt and Shmuel Onn.
\newblock On the diameter of convex polytopes.
\newblock {\em Discrete Mathematics}, 102:75--77, 1992.

\bibitem{MO}
Michal Melamed and Shmuel Onn.
\newblock Convex integer optimization by constantly many linear counterparts.
\newblock {\em Linear Algebra and its Applications}, 447:88--109, 2014.

\bibitem{N89}
Dennis Naddef.
\newblock The {H}irsch conjecture is true for $(0,1)$-polytopes.
\newblock {\em Mathematical Programming}, 45:109--110, 1989.

\bibitem{Onn}
Shmuel Onn.
\newblock {\em Nonlinear Discrete Optimization}.
\newblock Zurich Lectures in Advanced Mathematics. European Mathematical
  Society, 2010.

\bibitem{OR}
Shmuel Onn and Uriel~G. Rothblum.
\newblock Convex combinatorial optimization.
\newblock {\em Discrete and Computational Geometry}, 32:549--566, 2004.

\bibitem{S12}
Francisco Santos.
\newblock A counterexample to the {H}irsch conjecture.
\newblock {\em Annals of Mathematics}, 176:383--412, 2012.

\bibitem{OEI}
Neil {Sloane (editor)}.
\newblock The on-line encyclopedia of integer sequences.
\newblock \url{https://oeis.org}.

\bibitem{SS16}
Ivan Soprunov and Jenya Soprunova.
\newblock {Eventual quasi-linearity of the {M}inkowski length}.
\newblock {\em European Journal of Combinatorics}, 58:110--117, 2016.

\bibitem{S16}
Noriyoshi {Sukegawa}.
\newblock {Improving bounds on the diameter of a polyhedron in high
  dimensions}.
\newblock {\em arXiv:1604.04338}, 2016.

\bibitem{T91}
Torsten Thiele.
\newblock Extremalprobleme f\"{u}r {P}unktmengen.
\newblock {\em Diplomarbeit, Freie Universit\"{a}t Berlin}, 1991.

\bibitem{T14}
Michael {Todd}.
\newblock {An improved Kalai{-}Kleitman bound for the diameter of a
  polyhedron}.
\newblock {\em SIAM Journal on Discrete Mathematics}, 28:1944--1947, 2014.

\bibitem{Z95}
G{\"u}nter Ziegler.
\newblock {\em Lectures on Polytopes}.
\newblock Graduate Texts in Mathematics. Springer, 1995.

\end{thebibliography}

 $\:$\\
 \noindent {\small Antoine Deza,
{ Advanced Optimization Laboratory},  {Faculty of Engineering}\\
{ McMaster University}, {Hamilton, Ontario, Canada}. \\
{\em Email}: deza{\small @}mcmaster.ca}\\

 \noindent {\small George Manoussakis,
 Laboratoire de Recherche en Informatique\\
 Universit\'e de Paris Sud, Orsay, France.\\
{\em Email}: george{\small @}lri.fr}\\

 \noindent {\small Shmuel Onn,
 Operations Research, Davidson faculty of IE \& M\\
Technion Israel Institute of Technology, Haifa, Israel.\\
 {\em Email}: onn{\small @}ie.technion.ac.il}

\end{document}